
%


\hoffset -.5in

\font\bb=msym10

\font\twelvebf=cmbx12

\def\P{\hbox{\bb P}}
\def\E{\hbox{\bb E}}

\def\R{\hbox{\bb R}}

\openup 1\jot

\def\eps{\varepsilon}
\def\fl{\flushpar}
\def\pa{\partial}
\hfuzz=3cm
\documentstyle{amsppt}
\topmatter
\title
A Gordon-Chevet type inequality
\endtitle
\author 
B. Khaoulani   \footnote{Universit\'e Paris VII, URA 1321, 1990-91.}
\endauthor

\abstract We prove a new inequality for Gaussian processes, this
inequality implies the  Gordon-Chevet inequality. Some remarks on Gaussian
proofs of Dvoretzky's theorem are given.
\endabstract
\endtopmatter
\document
      
\fl\bf I. Introduction: \rm

Let $\{g_{i,k} \} (1 \leq i \leq n, 1 \leq k \leq d), \{h_k\}^d_1$ and
$\{g_i\}^n_1$
 denote independent sets of orthonormal Gaussian random variables.

Let $E$ and $F$ be Banach spaces,\  $\{f_k\}^d_{k=1} \subset F$ and
$\{x_i^*\}^n_{i=1} \subset E^*$.

Let $T(\omega) = \sum^n_{i=1} \sum^d_{k=1} g_{i,k} (\omega)
x_i^* \bigotimes f_k$ be
a random operator from $E$ to $F$. The Gordon-Chevet inequality says;
\cite{Cv}, \cite{G}.
$$
\align
\inf_{||x||_E=1}\{\left(\sum_{i=1}^nx^*_i(x)^2\right)^{\frac{1}{2}}\}
&\E (||\sum^d_{k=1} h_kf_k||)
- \varepsilon_2 (f_1,\cdots,f_d) \E (||\sum^n_{i=1}g_ix^*_i||_{E^*} \\
&\leq \E (\min_{||x||_E=1} ||T_{\omega}x||) \,\leq \,
\E (\max_{||x||_E=1} ||T_{\omega}x||)\\
&\leq \varepsilon_2 (x^*_1,\cdots,x^*_n)\E (||\sum^d_{k=1} h_k f_k||) +
\varepsilon_2 (f_1,\cdots,f_d)\E (||\sum^n_{i=1}g_ix^*_i||_{E^*})
\endalign
$$
where
$$
\varepsilon_2 (x^*_1,\cdots,x_n^*) = \sup\{ \left( \sum_{1\leq i\leq n}
x_i^*(x)^2\right)^{\frac{1}{2}} ; \Vert x \Vert_{E}\leq 1\}.
$$
and;
$$
\varepsilon_2 (f_1,\cdots,f_d) = \sup\{ \left( \sum_{1\leq k\leq d}y^*(f_k)^2
\right)^{\frac{1}{2}} ; \Vert y^*\Vert_{F^*}\leq 1\}.
$$

\fl The right-hand side inequality follows from  Chevet-inequality [Cv] and
can be
obtained from Sudakov lemma $ [G_1], [F]$ .

The left-hand side inequality is due to Gordon $[G_1]$ and follows from the
Gordon-Sudakov inequality.

Our aim is to deduce these inequalities from a general Gaussian inequality
for Gaussian processes.
\vskip 0.5cm
II.

Let $(\Omega,{\Cal A}, \P )$ be a probability space, $X$ be a canonical
$\R^d$-valued Gaussian random vector, (i.e. with covariance matrix
equal to $Id_d$). We define two Gaussian processes as follows.

For $n \geq 1$, let $B_2^n$ be the closed unit ball of $l^n_2$ and
 $S^{n-1}$ its unit sphere.
 For  $x=(x^1,\cdots,x^n)\in \R^n,$ let  $||x||_2=\left(\sum\limits_{i=1}^n\
(x^i)^2\right)^{\frac{1}{2}}$ and  let $X_1,\cdots,X_n$ be
$n$ independent copies of $X$,
independent
of $X$. Let $\{g_1,\cdots, g_n \}$ be a set of orthonormal
Gaussian random variables
independent of $ \{X,X_1,\cdots, X_n \}$.
$$
\text{Let} \,\,\, X_x = \sum^n_{i=1} x^i X_i \qquad \text{and} \qquad g_x =
\sum^n_{i=1}x^i g_i. \ \tag 2.1
$$
We shall prove the following inequality.

\fl \bf Theorem 1: \rm

Let $A \subset B_2^n$. Let $ F_x: \ \R^d \longrightarrow \R$ be a family of
1-Lipschitz functions indexed by  $x\in A$. Then the Gaussian
processes $\{ X_x\}_{x \in A},
\{g_x \}_{x \in A}$ satisfy
$$
\E \max_{x\in A} F_x(X_x) \leq \E \max_{x\in A} \{
F_x(||x||_2 X) + g_x \}. \ \tag 2.2
$$

\fl\bf Corollary 1: \rm

Let $A \subset B_2^n$, and $|||.|||$ be a norm on $\R^d$ such that;
 $\forall x \in \R^d \ |||x||| \leq ||x||_2,$
then the processes
$\{ X_x \}_{x \in A}$ and $\{g_x \}_{x \in A}$ verify

$$
\min_{x \in A}||x||_2\E |||X||| - \E\max_{x\in A}g_x \leq \E \min_{x \in A}
|||X_x||| \leq \E\max_{x \in A} |||X_x|||
\ \tag 2.3
$$

$$
\qquad \qquad \qquad \qquad \qquad \leq \E |||X||| + \E\max_{x\in A}g_x.
$$

\fl\bf Proof: \rm

For the right-hand side inequality, put: $F_y(x) = |||x|||$ and for
the left-hand
side inequality, put: $F_y(x) = -|||x|||$.

\fl\bf Corollary 2 \rm :
Let $X$ be a canonical $\R^d$-valued Gaussian random vector, $X_x$ and
$g_x$ are defined as in (2.1). Let  $A \subset S^{n-1}$,
$F$ be $1-$Lipschitz function on $\R^d$ and $\mu = \E F(X)$, then the
processes
$\{ X_x \}_{x \in A}$ and $\{g_x \}_{x \in A}$ verify :
$$
\E\max_{x\in A}|F(X_x)-\mu| \le \E |F(X)-\mu| +
\E\max_{x\in A}g_x \le 1 + \E\max_{x\in A}g_x.
$$
\fl \bf proof: \rm
For the first inequality, take $G(.)=|F(.)-\mu |$ which is a 1-Lipschitz
function, for the second, we use a well known Poincar\'e-type
inequality that is;

$$
\E |f(X)-\E (f(X))|^2 \le \E ||\bigtriangledown f(X)||_2^2,
$$
for $X$ as above, and all 1-Lipschitz function $f$ on $\R^d$,
\, \cite {P1},
\cite{C}.

\vskip .1in 
Next, we show how the Gordon-Chevet inequality follows from our inequality.
\vskip .1in 
\fl\bf Proof: \rm

Let $u:\R^d \rightarrow F$, $ u( \sum^d_{k=1} \alpha^k e_k ) = \sum^d_{k=1}
\alpha^k f_k$, and $ v: E \rightarrow l^2_n$ , $v(x) = (x^*_1(x),
\cdots ,x^*_n(x))$.

We have $ ||u|| = \varepsilon_2(f_1, \cdots, f_d),$ and $||v|| =
\varepsilon_2(x^*_1, \cdots, x^*_n).$

Let $X = \sum^d_{k=1} h_k e_k$ and for $ 1 \leq i \leq n$ let $ X_i = \sum^d_
{k=1} g_{ik} e_k$ then $X$ is an $\R^d$-valued canonical Gaussian vector and
$ X_1, \cdots, X_n \ n$ independent copies of $X$, independent of $X$.

Then; $\,u(X_{v(x)}(\omega)) = T_\omega(x)$, so the rest of the proof is as
in corollary 1 with $A=v(S_E),$  where $S_E$ is the unit sphere of $E$,
and $|||\alpha |||=||u(\alpha )||$.

\bigskip
Before proving theorem 1, we get a vectorial Slepian type inequality, from
which we deduce theorem 1, (see theorem 2 below).

\bigskip
We define some notations. For $x=(x_i)$, $y=(y_i)$ in  $\R^d$, $x\bigotimes y$
will denote
the matrix $\left(x_iy_j\right)_{1\le i,j\le d}$, and for $u,v \in \R^d$,
define $x\bigotimes y[u,v]$ as $<u,x\bigotimes y(v) > = <x,u><y,v>$, and
$||.||_{{\Cal L}(R^d)}$ the operator norm.
\bigskip

\bf Theorem 2 \rm

Let $\{ X_t\}$ , $\{ Y_t\}$, $t\in T$ be two families of Gaussian vectors with
values in $\R^d$, $\left\{g_t\right\}$ a family of Gaussian random variables
independent of $\{ X_t\}$ and  $\{ Y_t\}$, suppose

(i)\,\,\, $dist(X_t) = dist(Y_t)$ \,\,\, for all $t$ in $T,$

(ii)\,\,\, $||\E \left( X_t\bigotimes X_s - Y_t\bigotimes Y_s \right)
||_{{\Cal L}(\R^d)} \leq \frac{1}{2}\E|g_t - g_s|^2$, \, for all $s,t$ in $T$.

Let ${F_t}$, $t\in T$, be a family of  real 1-Lipschitz functions on $\R^d$,
 then;
 $$\E \sup_tF_t(X_t) \le \E \sup_t\{ F_t(Y_t) + g_t \}.$$

\fl\bf Proof \rm

We may clearly assume without loss of generality that the two processes
$\{X_t, t\in T\}$,
$\{Y_t, t\in T\}$ are independent, and also by a standard approximation
argument that the $F_t$ are 1-lipschitz and twice differentiable.

It is clear that we just need to prove the inequality for finite sets
${X_1,\cdots ,X_N}$, ${Y_1,\cdots ,Y_N}$, $(N\ge 1)$.

Fix  ${X_1,\cdots ,X_N}$ and  ${Y_1,\cdots ,Y_N}$, and prove that
$$
\E \max_{1 \leq i \leq N} \{ F_i(X_i) \} \leq \E \max_{1 \leq i \leq N} \{F_i
(Y_i)+ g_i \}. \ \tag 2.4
$$
For $ \theta \in [0, \frac{\pi}{2}]$ let
$$
Z(\theta) = (\cos (\theta )X_1 + \sin (\theta )Y_1,\sin (\theta ) g_1 ;
\cdots ; \cos (\theta ) X_N + \sin (\theta ) Y_N, \sin (\theta ) g_N)
$$
$Z(\theta ) $ is an $(\R^{d+1})^N$-valued Gaussian vector, with
$$
Z(0) = (X_1,0; \cdots; X_N,0) \qquad\text{and}\qquad
Z(\frac{\pi}{2}) = (Y_1,g_1; \cdots; Y_N,g_N);
$$
a vector $(y,z)$ of $E = (\R^{d+1})^N$ will be denoted by
$$
(y,z) = \left((y_i ; z_i)\right)_{1 \leq i \leq N} \,\text{where}\,\,\, y_i
\in \R^d \,\,\text{and}\,\, z_i \in \R.
$$
\vskip .1in

\fl \bf Step 1: \rm We prove the following lemma.
\vskip .1in 
\fl \bf Lemma: \rm

Let $F: \R^{(d+1)N} \rightarrow \R^N$, $F(y,z) = \left ( F_1(y_1)+z_1, \cdots ,
F_N(y_N)+z_N \right )$, \,where $F_1, \cdots ,F_N,$ are 1-Lipschitz on $\R^d,$
$G: \R^N \rightarrow \R$ be a twice differentiable function such that
$\exists k_1, k_2, $ s.t  $|G(.)| \le k_1e^{k_2||.||_2},$
$|\frac{\partial G(.)}{\partial \alpha_i}|\le k_1e^{k_2||.||_2}$ and,
$|\frac{\partial^2G(.)}{\partial \alpha_i\partial \alpha_j}|\le k_1e^{k_2
||.||_2}$ for all $i, j=1,\cdots ,N$.

Put; $\varphi = G\circ F$, and
$$
h(\theta ) = \E \varphi \left( Z(\theta ) \right). \ \tag 2.5
$$
Suppose;
$$
\forall i,j; \quad i \neq j \qquad \frac{\partial^2 G}{\partial
\alpha_i \partial
\alpha_j}
\leq 0, \ \tag 2.6
$$
$$
\text{and}\qquad\forall j = 1, \cdots,N \qquad \sum^N_{i=1} \frac{\partial^2
 G}{\partial\alpha_i \partial \alpha_j} = 0. \ \tag 2.7
$$
Then;

$h(\theta )$ is increasing, therefore
$$
\E G(F_1(X_1),\cdots ,F_N(X_N)) = h(0) \leq  h(\frac{\pi}{2}) =
\E G(F_1(Y_1) + g_1,\cdots ,F_N(Y_N) + g_N).
$$
\fl \bf Proof of the lemma. \rm
\vskip .1in

Let $\varepsilon > 0$, and $\Lambda$ an $(\R^{d+1})^N$-valued canonical
Gaussian vector independent of $\{ Z(\theta ); \theta \in ]0,
\frac{\pi}{2} [ \}$.

Let $Z_{\varepsilon} (\theta ) = Z(\theta ) + \varepsilon \Lambda$ so that
$\Gamma_{\varepsilon}(\theta ) = \Gamma(\theta ) + \varepsilon^2 I_E$  where
$\Gamma(\theta )$ is the covariance matrix of $Z(\theta )$ and $\Gamma_
{\varepsilon}(\theta )$ is the covariance matrix of $Z_{\varepsilon}(
\theta )$. Thus
$$
\Gamma_{\varepsilon}(\theta ) \longrightarrow_{\varepsilon \rightarrow 0}
\Gamma
(\theta ) \ \ \text{so that}\ \ h_{\varepsilon}(\theta )
\longrightarrow_{\varepsilon
 \rightarrow 0} h(\theta ).
$$
Remark  that
$$
\forall (u,v) \in E \quad < (u,v), \Gamma_{\varepsilon}(\theta ) (u,v) > \geq
\varepsilon^2 ||(u,v)||^2_E.
$$
Let $g_{\varepsilon} (y,z; \theta )$ be the density function of
$Z_{\varepsilon}
(\theta )$. We will list the following well-known identities:(see
$[G_2],[F],[G_1])$
$$
g_{\varepsilon} (y,z; \theta ) = \frac{1}{(2 \pi)^{(d+1)N}}
\int_E \exp\{i<(u,v);(y,z)> - \frac12 <(u,v),
\Gamma_{\varepsilon}(\theta )(u,v)>\}du dv
\tag 2.8
$$
where $du = du_1 \cdots du_N$, \ $du_i = du_{i,1} \cdots du_{i,d}$ and $dv =
dv_1 \cdots dv_N$
$$
h_{\varepsilon}(\theta ) = \int_E \varphi (y,z)
g_{\varepsilon}(y,z, \theta ) dy
dz \qquad \left( = \E \varphi(Z_{\varepsilon}(\theta ))\right) \ \tag 2.9
$$
$$
h{'}_{\varepsilon}(\theta ) = \int_E \varphi (y,z)
\frac{\partial}{\partial \theta }
g_{\varepsilon}(y,z, \theta ) dydz \ \tag 2.10
$$
$$
\frac{\partial}{\partial \theta}g_{\varepsilon}(x, \theta ) = \frac12 \sum_
{i,j=1}^{(d+1)N}\frac{d}{d \theta} \gamma_{i,j}^{\varepsilon}(\theta )
\frac{\partial^2}{\partial x_i \partial x_j} g_{\varepsilon}(x, \theta ) \
\tag 2.11
$$
where $x = (y,z) \quad$ and $\quad \Gamma_{\varepsilon}(\theta ) =
\left( \gamma_{i,j}^{
\varepsilon}(\theta ) \right)_{1 \leq i,j \leq N(d+1)}.$

\vskip .1in 

We compute $\Gamma_{\varepsilon} (\theta ).$

We can write $\Gamma_{\varepsilon}(\theta )$ as a block matrix : $\Gamma_{
\varepsilon}(\theta ) = \left(\Gamma_{i,j}^{\varepsilon}(\theta )\right)_{
1 \leq i \leq N , 1 \leq j \leq N }$ where
$$
\Gamma_{i,j}^{\varepsilon} (\theta ) = \E [Z_i^{\varepsilon}
(\theta )\bigotimes Z_j^{\varepsilon}(\theta) ]  \tag 2.12
$$
$$
Z_i^{\varepsilon}(\theta) = \left(
X_i(\theta) + Y_i(\theta)+\varepsilon \Lambda_i , g_i(\theta) + \varepsilon
\Lambda{'}_i \right)
$$
$$
\Lambda = ( \Lambda_i,\Lambda{'}_i)_{1 \leq i \leq N} \qquad \qquad
\Lambda_i = ( \Lambda^1_i, \cdots , \Lambda^d_i),
$$
$$
X_i(\theta)=\cos (\theta)X_i,\qquad Y_i(\theta)=\sin (\theta)Y_i, \qquad
\text{and}\qquad g_i(\theta)=\sin (\theta)g_i.
$$
Using the fact that $\{X_1, \cdots ,X_N\}$, $\{Y_1,\cdots , Y_N\}$ and
$\{g_1,
\cdots , g_N\}$ are independent processes, we find that;

$$
\Gamma_{i,j}^\varepsilon(\theta ) = \left[ \matrix A_{i,j}^\varepsilon
(\theta ) + \varepsilon^2 Id_d &0 \\0&B_{i,j}^\varepsilon(\theta )
\endmatrix \right] \tag 2.13
$$
where $A_{i,j}^\varepsilon(\theta)$ is a $d\times d$ matrix, and
$B_{i,j}^\varepsilon(\theta )$ is a scalar such that:

$$
A^{\varepsilon}_{ij}(\theta ) = \cos^2 (\theta )\E(X_i\bigotimes X_j) +
\sin^2 (\theta )\E(Y_i\bigotimes Y_j),
\quad\text{and}\quad B^{\varepsilon}_{ij}(\theta ) = \sin^2 (\theta )
\E g_ig_j + \varepsilon^2\delta_{i,j}.  \ \tag 2.14
$$
where $\delta_{i,j} = 1$ if $i=j$, and $0$ if $i\neq j.$

A simple computation gives
$$
\align
<(u,v); \Gamma_{\varepsilon} (\theta )(u,v)> &= \sum^N_{i=1} \sum^N_{j=1}
<u_i; A_{i,j}^{\varepsilon} (\theta ) u_j > \\
&+  \sum^N_{i=1} \sum^N_{j=1}
B_{i,j}^{\varepsilon}(\theta ) v_i.v_j.
\endalign
$$
Considering  $\frac{\partial^2}{\partial y_i \partial y_j}
g_{\varepsilon} (y,z; \theta)$ as a $d\times d$ matrix for each $i,j,$
this gives;
$$
\align
\frac{\partial}{\partial \theta} g_{\varepsilon}(y,z;\theta) &= \frac{1}{2}
\sum^N_{i,j=1}trace\left( \frac{\partial^2}{\partial y_i \partial y_j}
g_{\varepsilon} (y,z; \theta) \frac{d}{d \theta} A_{i,j}^{\varepsilon}
(\theta)\right) \\
&+ \frac{1}{2}  \sum^N_{i,j=1} \frac{d}{d \theta} B_{i,j}^{\varepsilon}
(\theta)\frac{\partial^2}{\partial z_i \partial z_j}
g_{\varepsilon} (y,z;\theta)
\endalign
$$

But
$$
h'_{\varepsilon}(\theta) = \int \varphi(y,z) \frac{\partial}{\partial \theta}
g_{\varepsilon}(y,z, \theta) dy dz  \tag 2.15
$$
Let  $M_{i,j}=\E Y_i\bigotimes Y_j - \E X_i\bigotimes X_j$, we get;

$$
h'_{\varepsilon}(\theta) = \frac{\sin 2 \theta}{2}\int_E \left \{
\sum^N_{i,j=1} trace\left(\frac{\partial^2 \varphi(y,z)}{\partial y_i
\partial y_j }.
M_{i,j}\right) + \sum^N_{i,j=1} \frac{\partial^2 \varphi(y,z)}
{\partial z_i \partial z_j}\E g_ig_j\right \} g_{\varepsilon}(y,z, \theta)
dy dz
 \tag 2.16
$$
Since \, $dist(X_i) = dist(Y_i)$ for all $i$, we get; $M_{i,i}=0$, hence we
have for $\varphi = G\circ F;$
$$
h'_{\varepsilon}(\theta) = \frac{\sin 2 \theta}{2}\int \left\{ \sum^N_
{i\neq j} tr\left(\frac{\pa^2 G\circ F}{\pa y_i \pa y_j}.M_{i,j}\right)  +
\sum^N_{i,j=1 }\left(\frac{\pa^2 G\circ F}{\pa z_i \pa z_j}\right) \E g_ig_j
\right\}g_{\varepsilon}(y,z;\theta) dy dz
$$
A simple computation gives for all $i\neq j$:
$$\frac{\pa^2G\circ F}{\pa y_i \pa y_j} = \frac{\pa^2G}{\pa \alpha_i
\pa \alpha_j}\circ F.\bigtriangledown F_i(y_i)\bigotimes \bigtriangledown
F_j(y_j)
$$
$$
\text{and}\qquad \frac{\pa^2G\circ F}{\pa z_i \pa z_j} = \frac{\pa^2G}
{\pa \alpha_i \pa \alpha_j}\circ F \,\qquad\text{for all}\,\, i,j
$$
Condition (2.7) gives;
$$
\frac{\partial^2 G}{\partial \alpha_i^2} = - \sum^N_{j=1
j \neq i} \frac{\partial^2 G}{\partial \alpha_i \partial \alpha_j} \,
\qquad\text{for all}
\,\, i,j
\ \tag 2.17
$$
so;
$
h'_{\varepsilon}(\theta)
$

$
= \frac{\sin 2 \theta}{2}\int \{ \sum^N_
{i\neq j} tr\left(\frac{\partial^2 G(F(y,z))}{\partial y_i \partial y_j}.
M_{i,j}\right) +  \sum^N_{i\neq j }\frac{\pa^2 G(F(y,z))}{\pa z_i \pa z_j}
\E g_ig_j$

$+\sum^N_{i=1}\frac{\partial^2G(F(y,z))}{\partial z_i^2}\E g_i^2 \}
\,\,\, g_{\varepsilon}(y,z;\theta) dy dz
$

$
= \frac{\sin 2 \theta}{2}\int \left\{ \sum^N_
{i\neq j} \left(tr\left(\frac{\pa^2 G(F(y,z))}{\pa y_i \pa y_j}.
M_{i,j}\right) + (\E g_ig_j - \frac{1}{2}[\E g_i^2+\E g_j^2])
\frac{\pa^2 G(F(y,z))}{\pa \alpha_i \pa \alpha_j} \right)\right\} 
g_{\varepsilon}(y,z;\theta) dy dz
$

$
=\frac{\sin 2 \theta}{2}\int \left\{ \sum^N_
{i\neq j} \left(\frac{\pa^2 G(F(y,z))}{\pa \alpha_i \pa \alpha_j}
\left(<M_{i,j}.\bigtriangledown F_i(y_i) , \bigtriangledown F_j(y_j)>
- \frac{1}{2}\E|g_i -g_j|^2 \right)\right)\right\} 
g_{\varepsilon}(y,z;\theta) dy  dz
$

Since  $||\vec{\bigtriangledown F}_i(y_i)|| \le 1$, then:

$<M_{i,j}({\bigtriangledown F}_i(y_i)), {\bigtriangledown F}_j(y_j)>
- \frac{1}{2}\E|g_i - g_j|^2$

$\le ||(M_{i,j})||_{{\Cal L}(\R^d)} - \frac{1}{2}\E|g_i - g_j|^2  \le 0.$

\vskip .1in 

so $h'_{\varepsilon }(\theta) \ge 0,$ and  $\E G(F(Z_{\varepsilon}(0)))
\le \E G(F(Z_{\varepsilon}(\frac{\pi}{2}))). \,$Finally, letting $\varepsilon
\to 0$, we get the result of the lemma.

\vskip .1in 
We finish the proof of theorem 2.

The map max satisfies the following conditions:

1) $\max (\alpha + t\vec 1) = \max(\alpha) +t \qquad \forall
\alpha \in \R^N$ and
$t \in \R$,  where $\vec 1 = (1,\cdots ,1) \in \R^N$.

2) it is slowly increasing and verifies in distribution sense \cite{G2}.
$$
\sum^N_{i=1} \frac{\partial \max}{\partial \alpha_i} = 1 \,\,, (2.6)\,
\text{and}\,\,\,(2.7)\; \tag 2.19
$$

So if we regularise max by convolution with a twice differentiable function
$\psi_k$ which is supported by a ball of radius $\frac{1}{k}$ we obtain a
function
$m_k$ which is 1-Lipschitz and satisfies the same conditions 1, 2, (2.6),
and (2.7). If k tends to infinity, we find by Lebesgue theorem that
$h$ is increasing in $[0; \frac{\pi}{2}]$, \,\, this completes  the proof.

\vskip 1,5cm

We prove theorem 1.

We have $X_x = \sum\limits_{i=1}^nx^iX_i$, let $Y_x= ||x||_2X$, where $x$
runs over a set $A\subset B_n^2$, then

$dist(X_x) = dist(Y_x)$.

Take a finite set $\{a_1,\cdots ,a_N\}$ in $A$, a simple computation gives ;
$$
M_{i,j} = \E\left( Y_{a_i}\bigotimes Y_{a_j} - X_{a_i}\bigotimes
X_{a_j}\right) =
\left( ||a_i||_2||a_j||_2 - a_i.a_j\right) Id_d,
$$
where $a_i.a_j$ is the scalar product.

Moreover, $\E |g_{a_i}-g_{a_j}|^2 = ||a_i - a_j||^2_2$, the $F_i$ are
1-Lipschitz functions so;

$
||M_{i,j}||_{{\Cal L}(\R^d)}- \frac{1}{2}\E|g_{a_i} - g_{a_j}|^2
$

$
= \left( ||a_i||_2||a_j||_2 - a_i.a_j\right) -\frac{1}{2}||a_i-a_j||_2^2
$

$
= -\frac{1}{2} \left( ||a_i||_2-||a_j||_2\right)^2 \le 0.
$

So conditions (i) and (ii) of theorem 1 are satisfied, and theorem 2 is
proved.

\vskip 1.5cm
III. Final remarks
\vskip .1in 
We give now a short proof of a   result due to V. Milman.
\line {}

\bf Theorem 3 \cite{M}, \cite{Sc} \rm.\,\,\, Let $\varepsilon >0$,
$f:\R^N \to \R$  be a lipschitz
function with constant L, $X=\sum\limits_{i=1}^Ng_ie_i$,
where $\{ g_i\}_{1\le i\le N}$ is a set of orthonormal
Gaussian random variables, $\{ e_i\}_{1\le i\le N}$ the canonical basis
of $l_2^N$, and $\mu = \E f(X)$. Then there exists an operator
$T:l^n_2\to \R^N$ \,with  $n = [\frac{\eps \mu}{L}
(\frac{\eps \mu}{L}-2)]$ , such that;
$$|f(Tx)-\mu| \le \eps \mu\,\,\text{for all}\,\, x\in S^{n-1}.$$

\vskip .1in 

\bf Proof \rm \,\,\, Consider  as above a r.v Gaussian operator $T_{\omega} =
\sum\limits_{i=1}^n \sum\limits_{j=1}^Ng_{ij}e_i^*\bigotimes e_j$
from $l^2_n$ to $\R^N$
$$X_i = \sum_{j=1}^Ng_{i,j}e_j, \,\,\,\text{and}\,\,\,
X_x=\sum_{i=1}^nx^iX_i.$$
Where $x=(x^i,\cdots ,x^n)$, then $X_x(\omega)=T_{\omega}x,$  we have;

$$
\align
\P\left(\left\{\omega / \exists x\in S^{n-1}; |f(X_x)-\mu | >
\eps\mu\right\}\right)
& = \P\left(\left\{ \omega ;\sup_{x\in S^{n-1}} |f(X_x)-\mu | >
\eps\mu\right\}\right) \\
& \le \frac{1}{\eps\mu}\E \sup_{x\in S^{n-1}} |f(X_x)-\mu |.
\endalign
$$
We apply corollary 2 to get;
$$
\P\left(\left\{\omega / \exists x\in S^{N-1}; |f(X_x)-\mu | >
\eps\mu\right\}\right)
\le \frac{1}{\eps\mu}\left\{ \E |f(X)-\mu | + L\E\sup_{x\in S^{n-1}}
\sum_{j=1}^nx^jg_j\right\}
$$
and using the Poincar\'e-type inequality as in corollary 2,
we find that;
$$
\P\left(\left\{ \omega / \exists x\in S^{n-1}; |f(X_x)-\mu | >
\eps\mu\right\}\right) \le \frac{L}{\eps\mu}[1+\E\sup_{x\in S^{n-1}}
\sum_{j=1}^nx^jg_j].
$$
$$ \qquad \quad \qquad \le \frac{L}{\varepsilon \mu }(1 + \sqrt{n}).
$$
We only need to choose $n$ such that this last expression is $< 1.$

\vskip 0.5cm

\vskip 0.5cm

Acknowlegment: \ I express my warmest thanks to Professor B. Maurey for
suggesting to me this question, and for many stimulating and fruitful
conversations.

\vskip 1cm

\fl \twelvebf References: \rm

\vskip 0.5cm

\itemitem [C]  \qquad L. Chen, A inequality for the multivariate normal
distributions, Journal Multivariate Anal. 12 (1982) 306-315.
\vskip 0.5cm
\itemitem [Cv] \qquad  S. Chevet, S\'eries
de variables al\'eatoires Gaussiennes
\'a valeurs
dans $E {\buildrel\wedge\over\bigotimes}_{\varepsilon} F$
Application aux produits d'espaces de Wiener,
S\'eminaire Maurey-Schwartz, expos\'e XiX  77/78.
\vskip 0.5cm
\itemitem [F] \qquad  X.M.Fernique, \ \  R\'egularit\'e
des trajectoires  des fonctions
al\'eatoires Gaussiennes, Lecture Notes in  Mathematics $n^0 480.$
\vskip 0.5cm
\itemitem [G 1] \qquad  Y.Gordon,\ \ Some inequalities for Gaussian
processes and
Applications. Israel J. Math. vol 50 (1985) p 265-289.
\vskip 0.5cm
\itemitem [G 2] \qquad  Y.Gordon, \ \ Elliptically contourned distributions
Prob. th. Rel. Field, 76, 429-438 (1987).
\vskip 0.5cm
\itemitem [K] \qquad J.P. Kahane,\ \ Une
in\'egalit\'e de type Slepian et Gordon
sur les processus Gaussiens, Israel J. Math. 55, p 109-110.
\vskip 0.5cm
\itemitem [M] \qquad V.D. Milman, New proof of the theorem of Dvoretzky on
sections of convex bodies, Funkcional Anal i Prilogen 5 (1971).
\vskip 0.5cm
\itemitem [P1] \qquad G. Pisier,\ \ Probabilistic methods in the geometry of
Banach spaces, Lecture Notes in Math; 1206  1-154 Springer, New York.
\vskip 0.5cm
\itemitem [Sc] \qquad G. Schechtman,\ \ A Remark concerning the dependence on
$\varepsilon$ in Dvoretzky's theorem. Springer LNM. 1376. p 274-277.

\end